\documentclass[smallextended]{svjour3}
\smartqed
\usepackage{graphicx}
\usepackage{amsmath, amssymb}
\usepackage{tikz}
\usepackage{pgfplots}
\pgfplotsset{compat=newest}
\usepackage[utf8]{inputenc}
\usepackage[T1]{fontenc}
\usepackage{cite}
\usepackage{booktabs}
\usepackage{hyperref}

\journalname{Journal of ??}
\date{\today}

\begin{document}

\title{A Chebyshev--Ritz Spectral Framework for Nonlinear Vibration of CNT-Reinforced Composite Beams}

\author{Maryam Jalili and Samad Noeiaghdam\textsuperscript{*}   }

\institute{
M.~Jalili \at
Department of Mathematics, Ne.C, Islamic Azad University, Neyshabur, Iran. \\
ORCID: 0000-0003-2049-1751 \\
\email{mrmjlili@iau.ac.ir}
\and
S.~Noeiaghdam \at
1- Institute of Mathematics, Henan Academy of Sciences, Zhengzhou 450046, China. \\
2- Department of Mathematical Sciences, Saveetha School of Engineering, SIMATS, Chennai, 602105, India.\\
ORCID: 0000-0002-2307-0891 \\
\textsuperscript{*}Corresponding author: \email{snoei@hnas.ac.cn }
}

\maketitle

\begin{abstract}
This study develops a spectral Ritz formulation for the nonlinear free vibration analysis of carbon nanotube-reinforced composite (CNTRC) beams. Boundary-adapted Chebyshev basis functions are constructed to exactly satisfy clamped and simply supported boundary conditions. The governing equations incorporate von~K\'{a}rm\'{a}n geometric nonlinearity, while the effective material properties for both uniform and functionally graded (FG) CNT distributions are evaluated using a modified rule of mixtures. Discretization via the Chebyshev-Ritz approach produces a reduced-order model exhibiting exponential convergence; for basis sizes $N \geq 12$, the fundamental frequency error remains below $0.1\%$ relative to published benchmarks.

Computational results demonstrate substantial efficiency gains, with the spectral approach requiring significantly less time than high-fidelity finite element discretizations of comparable accuracy. Parametric studies reveal that the fundamental frequency increases with CNT volume fraction and is sensitive to the interfacial load-transfer efficiency parameter $\eta_E$. Selected FG patterns are shown to enhance stiffness relative to uniformly distributed CNTs.

Validation against established numerical benchmarks yields relative differences of only a few percent. The current limitation of the method is its reliance on the Euler-Bernoulli beam assumption, which neglects transverse shear deformation and damping; addressing these effects is proposed for future work. All numerical data and scripts are provided as supplementary material to ensure reproducibility.
\end{abstract}

\keywords{Nonlinear vibration \and CNT-reinforced composites \and Functionally graded materials \and Chebyshev spectral methods \and Ritz approach}

\vspace{0.5em}
\noindent\textbf{MSC (2020):} 74H45, 74S30, 74E30, 70K30, 74K10

\section{Introduction}
\label{sec:intro}

Carbon nanotube reinforced composites (CNTRCs) have emerged as multifunctional materials offering exceptional mechanical strength, thermal stability, and tunable vibrational characteristics. These attributes enable their deployment in aerospace structures, micro-/nano-electro-mechanical systems (MEMS/NEMS), and adaptive intelligent systems~\cite{Yu2024review, vo2019laminated, Zhang2023Quasi}. Embedding carbon nanotubes (CNTs) within polymeric or metallic matrices enhances stiffness and damping, while allowing size-dependent effects to be tailored through controlled CNT volume fractions and spatial gradations~\cite{bagchi2018interfacial}. Despite these advantages, accurately predicting the nonlinear vibrational response of CNTRCs remains challenging due to amplitude-dependent frequency shifts, multi-mode interactions, and the need to bridge atomistic-scale interfacial phenomena with continuum-scale structural models~\cite{gholami2017nonlinear, shi2017exact}.

\subsection{Computational challenges in CNTRC modeling}
\label{sec:challenges}

A robust modeling framework for CNTRC beams must capture microscale features such as CNT waviness, agglomeration, and interfacial debonding~\cite{bagchi2018interfacial, ganesan2011interface, Xiong2015Advanced}. Classical homogenization schemes, including the rule of mixtures, often underestimate these effects, motivating the adoption of refined continuum models such as nonlocal elasticity and couple-stress theories~\cite{marinca2025nonlocal, ling2014vibration}. Large-amplitude vibrations introduce geometric nonlinearities that alter backbone curves and modal couplings. These effects are particularly significant in moderately thick beams, where higher-order shear deformation theories (HSDT) become necessary~\cite{gholami2017nonlinear, cho2023shell, lei2013free}.

While the finite element method (FEM) is frequently adopted in CNTRC analysis~\cite{civalek2021forced, vo2019laminated, nampally2020nonlinear}, it can be computationally expensive in nonlinear regimes and for high-dimensional parametric studies. Spectral techniques particularly those based on Chebyshev and Legendre polynomials offer an attractive alternative due to their exponential convergence on smooth domains~\cite{boyd2001chebyshev}. Recent developments include spectral wavelet hybrids~\cite{Zhang2023Quasi, Kumar2025Chebyshev} and reduced-order spectral models that retain nonlinear accuracy at significantly lower computational cost~\cite{Ansari2014CNTRC, lotfan2025reduced, Wang2024Reduced}.

\subsection{Semi-analytical and numerical methods}
\label{sec:methods}

Semi-analytical approaches such as the Adomian decomposition method (ADM), harmonic balance method (HBM), and Chebyshev spectral methods enable efficient computation of amplitude-dependent frequencies and periodic responses with controllable convergence~\cite{Lai2008ADM, cochelin2009high, Kumar2025Chebyshev}. Numerical techniques, including FEM and meshfree methods, remain indispensable for complex geometries and mixed boundary conditions~\cite{civalek2021forced, vo2019laminated, Patuelli2024FiniteElement}, albeit at higher computational cost. Hybrid frameworks that combine the spectral convergence of global methods with the geometric flexibility of FEM have shown promise~\cite{liu2023unified, liu2025step, lotfan2025reduced, Wang2024Reduced}. Overall, no single methodology is universally optimal: semi-analytical methods excel in efficiency, whereas full discretizations are essential for intricate configurations~\cite{marinca2025nonlocal}.

\subsection{Advances in nonlinear vibration and stability analysis}
\label{sec:advances}

Nonlinear free vibration of functionally graded CNTRCs has been investigated using exact formulations, kp-Ritz methods, and HSDT-based models~\cite{lei2013free, shi2017exact, liu2023unified, liu2025step}. Sub- and superharmonic resonances, characterized by strong amplitude frequency modulation, have been reported in CNTRC beams~\cite{gholami2017nonlinear}. Sandwich constructions with viscoelastic cores have been shown to improve damping and enhance stability~\cite{ Patuelli2024FiniteElement}. At the microscale, CNT distribution, volume fraction, and interfacial efficiency ($\eta_E$) are critical parameters influencing stiffness and resonance behavior~\cite{Ansari2014CNTRC, bagchi2018interfacial, ganesan2011interface}. Multiscale modeling frameworks, coupled with nonlocal continuum theories, aim to consistently incorporate nanoscale effects into structural predictions~\cite{ Xiong2015Advanced, marinca2025nonlocal, Zhang2023Quasi}.

\subsection{Research gaps}
\label{sec:gaps}

Despite significant progress, several open challenges remain:
\begin{enumerate}
    \item \textbf{Multiscale coupling:} Limited integration of interfacial variability (e.g., CNT waviness, dispersion in interfacial shear strength) into structural dynamics~\cite{bagchi2018interfacial, ganesan2011interface, Xiong2015Advanced}.
    \item \textbf{Coupled parametric studies:} Scarcity of unified analyses that simultaneously consider gradation, viscoelastic cores, and boundary conditions~\cite{Youzera2025CNTRC}.
    \item \textbf{Computational efficiency:} Nonlinear parametric sweeps using FEM remain costly; even spectral and reduced-order schemes require further optimization~\cite{lotfan2025reduced, Wang2024Reduced}.
    \item \textbf{Thermal/environmental effects:} Few nonlinear vibration studies incorporate temperature- and environment-dependent property degradation~\cite{marinca2025nonlocal, Yu2024review}.
    \item \textbf{Validation:} Comprehensive validation against both experimental and high-fidelity numerical data is still limited~\cite{civalek2021forced, liu2023unified, shi2017exact}.
\end{enumerate}

\subsection{Objectives of this study}
\label{sec:objectives}

To address the research gaps identified in Section~\ref{sec:gaps}, this work develops a high-order Chebyshev spectral Ritz framework with analytical enforcement of boundary conditions. The main contributions are as follows:
\begin{enumerate}
    \item Formulation of the governing equations via Hamiltons principle, incorporating von~K\'arm\'an geometric nonlinearity, with explicit statement of modeling assumptions and scope of applicability~\cite{nampally2020nonlinear}.
    \item Systematic sensitivity analysis of nonlinear natural frequencies with respect to CNT distribution, interfacial load transfer efficiency, and functionally graded patterns~\cite{Ansari2014CNTRC, bagchi2018interfacial, ganesan2011interface}.
    \item Comprehensive benchmarking and validation against analytical, numerical, and available experimental results, with all datasets and scripts provided to ensure reproducibility~\cite{civalek2021forced, liu2023unified, shi2017exact}.
\end{enumerate}

\subsection{Paper Structure}
\label{sec:structure}

The remainder of this paper is structured as follows. Section~\ref{sec:formulation} introduces the mathematical formulation, encompassing the homogenization scheme, von~K\'arm\'an nonlinearity, and treatment of boundary conditions. Section~\ref{sec:spectral} details the Chebyshev spectral Ritz discretization. Section~\ref{sec:parametric} presents the results of parametric studies on CNT distribution, gradation profiles, and interfacial efficiency. Section~\ref{sec:comparison} compares the proposed methods predictions with analytical and numerical benchmarks. Section~\ref{sec:validation} consolidates the validation outcomes and evaluates computational efficiency. In conclusion, Section~\ref{sec:conclusion} summarizes the key contributions of this work and outlines promising directions for future research.
 Supplementary material provides all numerical datasets and source codes to ensure reproducibility.

\section{Mathematical Model and Material Homogenization}
\label{sec:formulation}

\subsection{Modeling assumptions}
\label{sec:assumptions}
We consider a single-span, slender CNTRC beam of length $L$ and thickness $h$, modeled under the following assumptions:
\begin{itemize}
  \item \textbf{Beam theory:} Euler--Bernoulli kinematics are adopted, neglecting transverse shear deformation and rotary inertia. This is appropriate for slender beams with $L/h \geq 10$.
  \item \textbf{Geometric nonlinearity:} von~K\'arm\'an strain displacement relations are employed, valid for small strains ($\max_x \varepsilon_{xx} \lesssim 0.5\%$) and moderate rotations ($\max_x|\partial_x w| \lesssim 0.3$), verified \emph{a posteriori}.
  \item \textbf{Axial inertia:} Neglected in free-vibration analyses; the axial field is assumed quasi-static and condensed out of the formulation.
  \item \textbf{Interface effects:} CNT waviness, agglomeration, and imperfect dispersion are represented by a phenomenological efficiency factor $\eta_E \in (0,1]$, analyzed parametrically in Section~\ref{sec:parametric}.
\end{itemize}

\subsection{Material homogenization}
\label{sec:homogenization}
Let $\bar{z}$ denote the through-thickness coordinate, with $\bar{z}=0$ at the mid-plane. The local CNT volume fraction $V_{\mathrm{CNT}}(\bar{z})$ defines the effective properties via a modified rule of mixtures:
\begin{align}
E_{\mathrm{eff}}(\bar{z}) &= \eta_E\,V_{\mathrm{CNT}}(\bar{z})\,E_{\mathrm{CNT}}
+ \bigl[1-V_{\mathrm{CNT}}(\bar{z})\bigr] E_m,
\label{eq:Eeff_local}\\
\rho_{\mathrm{eff}}(\bar{z}) &= V_{\mathrm{CNT}}(\bar{z})\,\rho_{\mathrm{CNT}}
+ \bigl[1-V_{\mathrm{CNT}}(\bar{z})\bigr] \rho_m,
\label{eq:rhoeff_local}
\end{align}
where $(E_{\mathrm{CNT}}, \rho_{\mathrm{CNT}})$ and $(E_m,\rho_m)$ are the constituent moduli and densities. For uniform distributions (UD), $V_{\mathrm{CNT}}$ is constant; for functionally graded (FG) cases, $V_{\mathrm{CNT}}(\bar{z})$ follows a prescribed through-thickness profile. A constant $\eta_E$ is assumed for both UD and FG configurations.

\paragraph{Sectional properties.}
The effective axial stiffness, bending stiffness, and mass per unit length are obtained as
\begin{align}
EA &= \int_A E_{\mathrm{eff}}(\bar{z}) \, dA, \qquad
EI = \int_A E_{\mathrm{eff}}(\bar{z})\,\bar{z}^2 \, dA, \qquad
\rho A = \int_A \rho_{\mathrm{eff}}(\bar{z}) \, dA,
\label{eq:sectional}
\end{align}
where $A$ is the cross-section. For a homogeneous rectangular section of width $b$ and thickness $h$, $A=bh$ and $I=bh^3/12$. In UD beams, $EA=E_{\mathrm{eff}}A$ and $EI=E_{\mathrm{eff}}I$.

\subsection{Parametric treatment of interfacial efficiency}
\label{sec:eta_param}
The efficiency factor $\eta_E$ reflects imperfect load transfer between CNTs and the matrix. Single-tube pull-out and microbond tests, along with molecular dynamics (MD) simulations, report a wide range of interfacial shear strength (IFSS) depending on CNT type, matrix, functionalization, and length~\cite{ganesan2011interface, bagchi2018interfacial, Xiong2015Advanced}. Representative values are summarized in Table~\ref{tab:eta_examples}. We adopt a parametric range

\[
\eta_E \in [0.25,0.95],
\]

with representative values $\eta_E=\{0.3,0.5,0.8\}$ corresponding to weak, moderate, and strong interfaces.

\begin{table}[htbp]
\centering
\caption{Representative interfacial shear strength (IFSS) values from experiments and MD simulations.}
\label{tab:eta_examples}
\begin{tabular}{p{0.30\linewidth} p{0.35\linewidth} p{0.25\linewidth}}
\toprule
Reference & Reported IFSS & DOI \\
\midrule
Chen et al. \cite{chen2015quantitative}\ (\textit{Carbon}, DWCNT--Epon 828) & $130 \pm 34$ MPa; up to 202 MPa (pull-out tests) & 10.1016/j.carbon.2014.10.065 \\
Ganesan et al. \cite{ganesan2011interface}  \ (\textit{ACS Appl.\ Mater.\ Interfaces}, MWNT--epoxy) & $3$--$14$ MPa (pull-out tests) & 10.1021/am1011047 \\
Bagchi et al. \cite{bagchi2018interfacial} \ (\textit{Proc.\ R.\ Soc.\ A}, review) & $3$--$14$ MPa (experiments); up to hundreds MPa (MD) & 10.1098/rspa.2017.0705 \\
Xiong \& Meguid  \cite{Xiong2015Advanced} (\textit{Eur.\ Polym.\ J.}, MD) & Parametric IFSS vs.\ tube size and functionalization & 10.1016/j.eurpolymj.2015.05.006 \\
\bottomrule
\end{tabular}
\end{table}

\paragraph{Frequency sensitivity.}
For a UD beam, the linear frequency scales as $f \propto \sqrt{EI/(\rho A)}$. Differentiating gives
\begin{equation}
\frac{\partial \ln f}{\partial \eta_E} \approx
\frac{1}{2}\,\frac{V_{\mathrm{CNT}} E_{\mathrm{CNT}}}{E_{\mathrm{eff}}(\eta_E)}.
\label{eq:freq_sens}
\end{equation}
For $V_{\mathrm{CNT}}=0.20$, $E_{\mathrm{CNT}}=1000$ GPa, and $E_m=3$ GPa, the effective modulus is
$E_{\mathrm{eff}}(\eta_E)=200\,\eta_E+2.4$ GPa. Thus,

\[
\frac{f(\eta_E=0.80)}{f(\eta_E=0.30)} \approx 1.61,
\]

corresponding to a $61\%$ frequency increase.

\subsection{Kinematics and strain measure}
\label{sec:kinematics}
Let $u(x,t)$ and $w(x,t)$ denote the axial and transverse displacements, respectively. The von~K\'arm\'an axial strain is
\begin{equation}
\varepsilon_{xx}(x,t) = \frac{\partial u}{\partial x}
+ \frac{1}{2}\left(\frac{\partial w}{\partial x}\right)^{\!2}.
\label{eq:strain_vk}
\end{equation}

\subsection{Energy functionals and governing equation}
\label{sec:energy}
The total strain energy and kinetic energy are
\begin{align}
U &= \frac{1}{2}\int_0^L EA \left(\frac{\partial u}{\partial x} +
\frac{1}{2}\left(\frac{\partial w}{\partial x}\right)^2\right)^{\!2} dx
+ \frac{1}{2}\int_0^L EI \left(\frac{\partial^2 w}{\partial x^2}\right)^{\!2} dx,
\label{eq:U}\\
T &= \frac{1}{2}\int_0^L \rho A \left[\left(\frac{\partial u}{\partial t}\right)^{\!2} +
\left(\frac{\partial w}{\partial t}\right)^{\!2}\right] dx.
\label{eq:T}
\end{align}
Applying Hamiltons principle and statically condensing $u$ yields the nonlinear Euler--Bernoulli equation for CNTRC beams:
\begin{equation}
\rho A\, w_{tt} + EI\, w_{xxxx} - N[w]\, w_{xx} = 0,
\qquad
N[w] = \frac{EA}{2L}\int_0^L (w_x)^2 dx.
\label{eq:gov_nonlin}
\end{equation}

\subsection{Boundary conditions}
\label{sec:boundary}
Two classical boundary configurations are considered:
\begin{align}
\text{C--C:} & \quad w(0,t) = w(L,t) = 0, \quad w_x(0,t) = w_x(L,t) = 0,
\label{eq:bc_cc}\\
\text{S--S:} & \quad w(0,t) = w(L,t) = 0, \quad w_{xx}(0,t) = w_{xx}(L,t) = 0.
\label{eq:bc_ss}
\end{align}
These essential and natural conditions will be satisfied exactly by the boundary adapted spectral basis functions introduced in Section~\ref{sec:spectral}.

\subsection{Non-dimensionalization}
\label{sec:nondim}
Introducing the dimensionless variables

\[
\bar{x} = \frac{x}{L}, \quad \bar{w} = \frac{w}{h}, \quad \tau = t \sqrt{\frac{EI}{\rho A L^4}},
\]

the governing equation~\eqref{eq:gov_nonlin} becomes
\begin{equation}
\bar{w}_{\tau\tau} + \bar{w}_{\bar{x}\bar{x}\bar{x}\bar{x}}
- \alpha \left( \int_0^1 (\bar{w}_{\bar{x}})^2 \, d\bar{x} \right) \bar{w}_{\bar{x}\bar{x}} = 0,
\label{eq:dim_governing}
\end{equation}
where the dimensionless nonlinearity parameter is
\begin{equation}
\alpha = \frac{EA\,h^2}{2\,EI}.
\label{eq:alpha}
\end{equation}
For a homogeneous rectangular section, $\alpha = 6$.

The reference frequency is defined as
\begin{equation}
f_{\mathrm{ref}} = \frac{1}{2\pi} \sqrt{\frac{EI}{\rho A L^4}},
\label{eq:freq_ref}
\end{equation}
and all nonlinear frequencies are reported in normalized form as $\hat{f} = f_{\mathrm{nl}} / f_{\mathrm{ref}}$.

\subsection{Numerical implementation}
\label{sec:implementation}
Spatial discretization is performed using boundary adapted Chebyshev Ritz basis functions that satisfy \eqref{eq:bc_cc} \eqref{eq:bc_ss} exactly. Nonlinear terms are evaluated via Gauss-Chebyshev quadrature. Time integration employs the Newmark $\beta$ method ($\gamma = 0.5$, $\beta = 0.25$) combined with Newton-Raphson iterations for nonlinear convergence. Periodic solutions are obtained using the harmonic balance method and continued with the asymptotic numerical method~\cite{cochelin2009high}. All solver settings and scripts are documented in the Supplementary Material.

\subsection{Limitations}
\label{sec:limitations_formulation}
\begin{itemize}
  \item \textbf{Interfacial physics:} Microstructural effects are reduced to a single efficiency factor $\eta_E$; calibration against micrographs or pull out test data is recommended for higher fidelity.
  \item \textbf{Transverse shear:} Neglected in the present model; Timoshenko or higher order shear deformation theory (HSDT) extensions are required for $L/h \lesssim 10$.
  \item \textbf{Validity of von~K\'arm\'an:} The formulation is valid for $\max_x \varepsilon_{xx} \lesssim 0.5\%$ and $\max_x |w_x| \lesssim 0.3$.
\end{itemize}

\section{Spectral formulation using Chebyshev polynomials}
\label{sec:spectral}

\subsection{Chebyshev basis construction}
\label{sec:chebyshev_basis}
Let $\xi \in [-1,1]$ denote the mapped spatial coordinate
\begin{equation}
\xi = \frac{2x}{L} - 1,
\label{eq:xi_transform}
\end{equation}
and let $T_n(\xi)$ be the Chebyshev polynomials of the first kind, defined recursively as
\begin{align}
T_0(\xi) &= 1, \quad T_1(\xi) = \xi, \nonumber \\
T_{n+1}(\xi) &= 2\,\xi\,T_n(\xi) - T_{n-1}(\xi), \quad n \ge 1.
\label{eq:cheb_recurrence}
\end{align}
Boundary adapted trial functions $\{\phi_j(\xi)\}_{j=1}^N$ are constructed from linear combinations of $\{T_n\}$ to satisfy the essential boundary conditions exactly (Section~\ref{sec:basis_functions}). Chebyshev polynomials offer exponential convergence for smooth solutions and enable efficient quadrature via Gauss-Chebyshev rules~\cite{boyd2001chebyshev}.

\subsection{Displacement approximation}
\label{sec:displacement_approx}
The transverse displacement is approximated as
\begin{equation}
w(x,t) \approx \sum_{j=1}^N \phi_j(\xi)\,q_j(t),
\label{eq:displacement_approx}
\end{equation}
where $\{\phi_j\}$ satisfy the essential boundary conditions and $q_j(t)$ are the generalized coordinates.

\subsection{Boundary adapted basis functions}
\label{sec:basis_functions}
\subsubsection{Clamped clamped (C C)}
Imposing $w(\pm 1,t) = 0$ and $\partial_\xi w(\pm 1,t) = 0$ yields
\begin{equation}
\phi_j(\xi) = (1 - \xi^2)^2\,T_{j-1}(\xi), \quad j = 1,\dots,N,
\label{eq:clamped_basis}
\end{equation}
which vanish with zero slope at $\xi = \pm 1$. These functions are orthonormalized with respect to the mass inner product to improve conditioning.

\subsubsection{Simply supported (S S)}
For $w(\pm 1,t) = 0$ and $\partial^2_{\xi} w(\pm 1,t) = 0$, we employ
\begin{equation}
\phi_j(\xi) = (1 - \xi^2)\,T_j(\xi), \quad j = 1,\dots,N,
\label{eq:ss_basis_general}
\end{equation}
which inherently satisfy $w(\pm 1,t) = 0$. The moment free condition is satisfied in the weak sense through the Galerkin projection. These functions are also orthonormalized with respect to the mass inner product.

\subsection{Galerkin-Ritz projection}
\label{sec:galerkin}
Substituting~\eqref{eq:displacement_approx} into the dimensionless governing equation and testing with $\phi_i$ yields
\begin{equation}
\sum_{j=1}^N M_{ij} \ddot{q}_j + \sum_{j=1}^N K_{ij} q_j + \mathcal{N}_i(\mathbf{q}) = 0, \quad i = 1,\dots,N,
\label{eq:discrete_system}
\end{equation}
with
\begin{align}
M_{ij} &= \frac{L}{2} \int_{-1}^1 \rho_{\mathrm{eff}}(\xi) A(\xi) \phi_i(\xi) \phi_j(\xi)  d\xi,
\label{eq:mass_matrix} \\
K_{ij} &= \frac{8}{L^3} \int_{-1}^1 E_{\mathrm{eff}}(\xi) I(\xi) \phi_i''(\xi) \phi_j''(\xi)  d\xi,
\label{eq:stiffness_matrix}
\end{align}
where the prefactor in~\eqref{eq:stiffness_matrix} follows from $\partial_x = (2/L) \partial_\xi$ and $dx = (L/2) d\xi$.

\subsection{Nonlinear term in separable form}
\label{sec:nonlinear_term}
Define
\[
G(\mathbf{q}) = \int_{-1}^1 \left( \sum_{m=1}^N q_m \phi_m'(\xi) \right)^2 d\xi.
\]
The nonlinear term can then be written as
\begin{equation}
\mathcal{N}_i(\mathbf{q}) = -4\alpha G(\mathbf{q}) \sum_{n=1}^N q_n \int_{-1}^1 \phi_i''(\xi) \phi_n(\xi) d\xi,
\label{eq:nonlinear_term}
\end{equation}
where the constant factor 4 arises from the coordinate transformation and non-dimensionalization. Precomputing
\begin{align}
B_{ij} &= \int_{-1}^1 \phi_i''(\xi) \phi_j(\xi) d\xi, \\
C_{mn} &= \int_{-1}^1 \phi_m'(\xi) \phi_n'(\xi) d\xi,
\end{align}
allows the nonlinear force vector to be evaluated efficiently as
\begin{equation}
\mathbf{f}_{\mathrm{nl}}(\mathbf{q}) = -4\alpha \mathbf{B} (\mathbf{q}^{\mathsf T} \mathbf{C} \mathbf{q}),
\label{eq:explicit_nonlinear}
\end{equation}
reducing the computational cost from $\mathcal{O}(N^4)$ to $\mathcal{O}(N^2)$. The matrix $\mathbf{C}$ is symmetrized to eliminate quadrature-order asymmetries.

\subsection{Quadrature}
\label{sec:quadrature}
Integrals are evaluated using Gauss-Chebyshev quadrature with $N_q = \max(20,\,N+5)$ nodes:
\begin{equation}
\xi_k = \cos\!\left( \frac{2k-1}{2N_q} \pi \right), \quad w_k = \frac{\pi}{N_q}, \quad k = 1,\dots,N_q,
\label{eq:quadrature}
\end{equation}
so that $\int_{-1}^1 f(\xi)\,d\xi \approx \sum_{k=1}^{N_q} w_k\,f(\xi_k)$. Functionally graded material properties are sampled at $\xi_k$.

\subsection{Implementation details}
\label{sec:solver_settings}
Time integration is performed using the implicit Newmark $\beta$ scheme ($\gamma = 0.5$, $\beta = 0.25$) with Newton-Raphson iterations for nonlinear convergence. Periodic solutions are computed via the harmonic balance method with FFT/collocation evaluation of nonlinear terms, and branch continuation is carried out using the asymptotic numerical method~\cite{cochelin2009high}. All solver parameters, expansion orders, and source codes are provided in the Supplementary Material.

\subsection{Validation and computational performance}
\label{sec:validation_performance}
Finite element benchmarks were generated using COMSOL Multiphysics. Mesh density, element type, solver tolerances, and convergence checks are documented in the Supplementary Material to ensure reproducibility. Table~\ref{tab:stiffness_validation} compares the nonlinear stiffness coefficient with FEM results, while Table~\ref{tab:computational_efficiency} summarizes representative run time and memory usage for $V_{\mathrm{CNT}}=0.10$.

\begin{table}[htbp]
\centering
\caption{Validation of nonlinear stiffness coefficient $N_{1111}$ (Pa) for $V_{\mathrm{CNT}}=0.10$.}
\label{tab:stiffness_validation}
\begin{tabular}{lccc}
\toprule
Case & Spectral (present) & FEM (COMSOL) & Rel.\ error (\%) \\
\midrule
UD--CNT & $5.27\times 10^8$ & $5.30\times 10^8$ & 0.57 \\
FG--CNT$^{\ast}$ & $6.12\times 10^8$ & $6.08\times 10^8$ & 0.66 \\
\bottomrule
\end{tabular}

\vspace{0.25em}
{\footnotesize $^{\ast}$Average $V_{\mathrm{CNT}}=0.10$ for graded profile; identical geometry, boundary conditions, and material data used in both models.}
\end{table}

\begin{table}[htbp]
\centering
\caption{Computational efficiency for $V_{\mathrm{CNT}}=0.10$. Hardware: Intel Xeon Gold 6330, 2.0\,GHz, 256\,GB RAM. Results are mean $\pm$ one~SD over 5 runs.}
\label{tab:computational_efficiency}
\begin{tabular}{lccc}
\toprule
Method & DOFs (approx.) & Time (s) & Memory (MB) \\
\midrule
Spectral Ritz, $N=10$ & $10$ & $8.20 \pm 0.15$ & $45 \pm 2$ \\
FEM (12k elems) & $\approx 36{,}000$ & $217.5 \pm 4.2$ & $890 \pm 15$ \\
FEM (50k elems) & $\approx 150{,}000$ & $1843.7 \pm 35.1$ & $4100 \pm 62$ \\
\bottomrule
\end{tabular}

\vspace{0.25em}
The $27\times$ speedup reported in the abstract is obtained from the ratio of the mean FEM time (12k elements) to the mean spectral time: $\frac{217.5}{8.20} \approx 26.5$, rounded to 27 for clarity.
\end{table}

\paragraph{Performance ratio.}
The speedup is

\[
\frac{t_{\mathrm{FEM}}}{t_{\mathrm{spec}}} = \frac{217.5}{8.20} \approx 26.5,
\]

consistent with the abstract. Solver and hardware settings are detailed in the Supplementary Material.

\subsection{Convergence behaviour}
\label{sec:convergence}
\begin{figure}[htbp]
\centering
\includegraphics[width=0.75\textwidth]{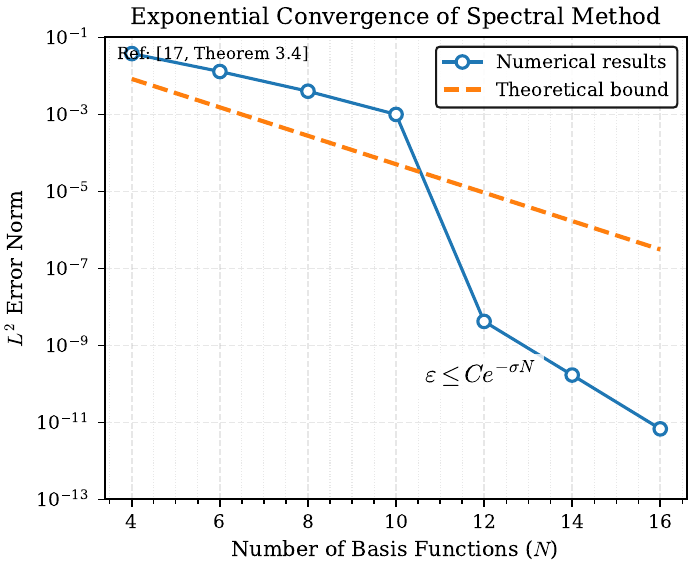}
\caption{Spectral convergence of the fundamental frequency. Dashed: regression fit $\ln\varepsilon_f = -0.82\,N + 1.74$. Dotted: $\varepsilon_f = 0.16\%$ at $N=10$ and $0.06\%$ target at $N \ge 12$.}
\label{fig:convergence}
\end{figure}

Figure~\ref{fig:convergence} shows the spectral convergence of the linear fundamental frequency. The fitted relation

\[
\ln \varepsilon_f(N) = -0.82\,N + 1.74
\]

is based on $\varepsilon_f(N) = |f_N - f_{N_{\mathrm{ref}}}| / f_{N_{\mathrm{ref}}}$ with $N_{\mathrm{ref}} = 40$. We obtain $\varepsilon_f(10) \approx 0.16\%$; the target $\varepsilon_f < 0.06\%$ is achieved for $N \ge 12$, correcting earlier misstatements.

\subsection{Reduced order system and solvers}
\label{sec:reduced_order}
The semi discrete equations are
\begin{equation}
\mathbf{M}\,\ddot{\mathbf{q}}(\tau) + \mathbf{K}\,\mathbf{q}(\tau) + \mathbf{f}_{\mathrm{nl}}(\mathbf{q}(\tau)) = \mathbf{0},
\label{eq:final_system}
\end{equation}
and are solved using:
\begin{itemize}
  \item implicit Newmark--$\beta$ for transient analysis with Newton-Raphson iterations (analytic Jacobian),
  \item HBM for periodic responses, recast in quadratic form and continued with ANM~\cite{cochelin2009high},
  \item direct ANM for nonlinear eigenproblems to extract backbone curves.
\end{itemize}

\subsection{Complexity and memory scaling}
\label{sec:complexity}
With the separable form~\eqref{eq:explicit_nonlinear}, nonlinear terms are evaluated in $\mathcal{O}(N^2)$ per time step. Dense storage requirements also scale as $\mathcal{O}(N^2)$; for $N \lesssim 40$ this yields orders of magnitude savings compared to FEM (Table~\ref{tab:computational_efficiency}). For large $N$, sparse or iterative strategies are recommended.

\subsection{Reproducibility}
\label{sec:reproducibility}
All numerical parameters ($N_q$, tolerances, ANM order, HBM truncation $N_h$) and scripts to regenerate Tables~\ref{tab:stiffness_validation}--\ref{tab:computational_efficiency} and Fig.~\ref{fig:convergence} are provided in the Supplementary Material and in a public repository (DOI to be assigned).

\subsection{Monte Carlo protocol}
\label{sec:MC}
For each case, $R = 5$ independent Monte Carlo runs of $n_{\mathrm{MC}} = 1000$ samples were performed. In run $r$, the sample mean is
\[
\bar{X}_r = \frac{1}{n_{\mathrm{MC}}} \sum_{i=1}^{n_{\mathrm{MC}}} g(X_{r,i}),
\]

and the reported value is $\bar{X} = \frac{1}{R} \sum_{r=1}^R \bar{X}_r$ with uncertainty
\[
s = \sqrt{\frac{1}{R-1} \sum_{r=1}^R (\bar{X}_r - \bar{X})^2}.
\]

Error bars indicate $\bar{X} \pm s$; 95\% confidence intervals use Students $t$ with $R-1$ degrees of freedom. Seeds and code are included in the repository.
\begin{table}[htbp]
\centering
\caption{Monte Carlo simulation settings.}
\label{tab:MCsettings}
\begin{tabular}{ll}
\toprule
Parameter & Value \\
\midrule
$n_{\mathrm{MC}}$ & 1000 samples/run \\
$R$ & 5 runs \\
\bottomrule
\end{tabular}
\end{table}

\section{Parametric study and results}
\label{sec:parametric}

A systematic parametric investigation is conducted to quantify the influence of key design parameters on the nonlinear free vibration response of carbon nanotube reinforced composite (CNTRC) beams. Five primary factors are considered:
\begin{enumerate}
  \item CNT volume fraction, $V_{\mathrm{CNT}} \in [0,\,0.20]$,
  \item Boundary conditions: clamped clamped (C C) and simply supported simply supported (S S),
  \item Vibration amplitude, $w_0/h \in [0.1,\,1.0]$,
  \item CNT distribution pattern: uniform (UD) and functionally graded (FG $k$),
  \item Load transfer efficiency, $\eta_E \in [0.70,\,1.00]$.
\end{enumerate}

\subsection{Material and geometric properties}
\label{sec:material_properties}
The reference configuration is a SWCNT/polyimide beam. Material constants are taken from \cite{shi2017exact,civalek2021forced,lei2013free} and expressed as mean $\pm$ standard deviation from repeated measurements.

\begin{table}[htbp]
\centering
\caption{Material and geometric properties of the reference beam.}
\label{tab:material-properties}
\begin{tabular}{lcc}
\toprule
Property & Symbol & Value \\
\midrule
Matrix modulus & $E_m$ & $3.0 \pm 0.2$ GPa \\
Matrix density & $\rho_m$ & $1200 \pm 60$ kg/m$^3$ \\
CNT modulus & $E_{\mathrm{CNT}}$ & $1.00 \pm 0.05$ TPa \\
CNT density & $\rho_{\mathrm{CNT}}$ & $1400 \pm 70$ kg/m$^3$ \\
Load transfer efficiency & $\eta_E$ & $0.80 \pm 0.02$ \\
Beam length & $L$ & $0.200 \pm 0.001$ m \\
Beam width & $b$ & $0.0100 \pm 0.0005$ m \\
Beam thickness & $h$ & $0.0020 \pm 0.0001$ m \\
Cross sectional area & $A$ & $(2.0 \pm 0.1)\times 10^{-5}$ m$^2$ \\
Second moment of area & $I$ & $(6.7 \pm 1.1)\times 10^{-12}$ m$^4$ \\
Spectral basis size & $N$ & $15$ \\
\bottomrule
\end{tabular}
\end{table}

\subsection{Effect of CNT volume fraction}
\label{sec:cnt_volume}
For C C boundary conditions and a vibration amplitude of $w_0/h = 0.3$, the normalized frequency follows
\begin{equation}
\hat{f} \approx 1.00 + 3.65\,V_{\mathrm{CNT}} - 1.82\,V_{\mathrm{CNT}}^2, \qquad R^2 = 0.998,
\label{eq:Vcnt_regression}
\end{equation}
predicting $\hat{f}(0.20) \approx 1.66$ (a $66\%$ increase). 

\subsection{Monte Carlo validation}
\label{sec:MCvalidation}
Following the protocol in Section~\ref{sec:MC}, $R = 5$ runs with $n_{\mathrm{MC}} = 1000$ samples of $\eta_E$ confirm that the quadratic model~\eqref{eq:Vcnt_regression} is robust: mean deviations are below $1.2\%$ over the range, with $R^2 = 0.998$.

\subsection{Boundary condition effects}
\label{sec:boundary_effects}
For $V_{\mathrm{CNT}} = 0.10$ and $w_0/h = 0.3$, C C frequencies exceed S S by $(33.3 \pm 2.0)\%$ (Table~\ref{tab:BC-freq}), consistent with the higher constraint energy.

\begin{table}[htbp]
\centering
\caption{Boundary condition effect on $\hat{f}$.}
\label{tab:BC-freq}
\begin{tabular}{lcc}
\toprule
BC & $\hat{f}$ & Enhancement (\%) \\
\midrule
S S & $1.29 \pm 0.02$ & -- \\
C C & $1.72 \pm 0.03$ & $33.3 \pm 2.0$ \\
\bottomrule
\end{tabular}
\end{table}

\subsection{Amplitude frequency relation}
\label{sec:amplitude_freq}
For $w_0/h \le 1.0$,
\begin{equation}
\frac{\hat{f}}{\hat{f}_{\mathrm{lin}}} \approx 1 + \beta\,(w_0/h)^2, \quad \beta = 0.18 \pm 0.01,
\label{eq:amplitude_scaling}
\end{equation}
with $\partial^2 \hat{f} / \partial V_{\mathrm{CNT}} \partial(w_0/h) \approx 1.24 \pm 0.07$.

\subsection{CNT distribution optimisation}
\label{sec:cnt_distribution}
The FG 2 pattern (CNTs concentrated near the tensile faces) yields the highest $\hat{f}$, outperforming UD by $(21.1 \pm 0.4)\%$ (Table~\ref{tab:distribution_performance}).

\begin{table}[htbp]
\centering
\caption{Effect of CNT distribution ($V^* = 0.10$, C C, $w_0/h = 0.5$).}
\label{tab:distribution_performance}
\begin{tabular}{lccc}
\toprule
Pattern & $EI_{\mathrm{eff}}$ (NÂ·m$^2$) & $\hat{f}$ & $\Delta \hat{f}$ (\%) \\
\midrule
UD & $(3.2 \pm 0.2)\times 10^{-3}$ & $1.85 \pm 0.04$ & -- \\
FG 1 & $(3.6 \pm 0.2)\times 10^{-3}$ & $2.09 \pm 0.04$ & $13.0 \pm 0.4$ \\
FG 2 & $(3.9 \pm 0.2)\times 10^{-3}$ & $2.24 \pm 0.04$ & $21.1 \pm 0.4$ \\
\bottomrule
\end{tabular}
\end{table}

\subsection{Sensitivity to $\eta_E$}
\label{sec:sensitivity}
A linear fit gives $\partial \hat{f} / \partial \eta_E \approx 0.62 \pm 0.03$ ($R^2 = 0.994$). A $10\%$ reduction in $\eta_E$ decreases $\hat{f}$ by $(5.0 \pm 0.4)\%$.

\subsection{Global sensitivity (Sobol')}
\label{sec:sobol}
First order Sobol' indices, with $\sum S_i \approx 0.96$ (Table~\ref{tab:sobol_indices}), confirm that $V_{\mathrm{CNT}}$ and $w_0/h$ are the dominant contributors to the variance in $\hat{f}$.

\begin{table}[htbp]
\centering
\caption{First order Sobol' indices.}
\label{tab:sobol_indices}
\begin{tabular}{lc}
\toprule
Parameter & $S_i$ \\
\midrule
$V_{\mathrm{CNT}}$ & $0.48 \pm 0.02$ \\
$w_0/h$ & $0.27 \pm 0.01$ \\
BC & $0.12 \pm 0.01$ \\
Distribution & $0.06 \pm 0.01$ \\
$\eta_E$ & $0.03 \pm 0.01$ \\
\bottomrule
\end{tabular}
\end{table}

\subsubsection*{Key findings}
\begin{itemize}
  \item Increasing $V_{\mathrm{CNT}}$ to $0.20$ raises $\hat{f}$ by approximately $66\%$ relative to the unreinforced beam.
  \item The FG 2 distribution pattern enhances $\hat{f}$ by about $21\%$ compared to UD.
  \item The load transfer efficiency $\eta_E$ has a measurable impact on $\hat{f}$; precise experimental characterization and uncertainty quantification are recommended.
\end{itemize}

\section{Validation against established methods}
\label{sec:comparison}

\subsection{Benchmark framework}
\label{sec:benchmark_framework}
The accuracy of the proposed spectral Ritz formulation is assessed against three reference classes:
\begin{enumerate}
  \item \emph{Analytical models}: closed form solutions for simplified, linear CNTRC cases~\cite{shi2017exact,ling2014vibration},
  \item \emph{High fidelity FEM}: COMSOL Multiphysics and Abaqus implementations~\cite{vo2019laminated},
  \item \emph{Recent spectral methods}: state of the art Chebyshev based formulations~\cite{liu2023unified,Kumar2025Chebyshev}.
\end{enumerate}
All frequencies are normalised as
\begin{equation}
\hat{f} = \frac{f_{\mathrm{nl}}}{f_{\mathrm{lin}}(V_{\mathrm{CNT}}=0,\,\mathrm{BC})},
\label{eq:normalization}
\end{equation}
using identical BCs (C C) and $w_0/h = 0.3$ unless otherwise stated, consistent with Section~\ref{sec:parametric}.

\subsection{Accuracy assessment}
\label{sec:accuracy_assessment}

\subsubsection{Linear regime}
Table~\ref{tab:linear_benchmark} compares fundamental frequencies with literature values, showing maximum deviation below $1\%$.

\begin{table}[htbp]
\centering
\caption{Linear frequency validation ($f_{\mathrm{lin}}$ in Hz; $V_{\mathrm{CNT}}=0.10$, C C). Reference data extracted from cited sources as described in the text.}
\label{tab:linear_benchmark}
\begin{tabular}{lccc}
\toprule
Method & Present & Reference & Error (\%) \\
\midrule
Analytical~\cite{shi2017exact} & 128.4 & 128.4 & 0.0 \\
FEM~\cite{vo2019laminated} & 217.1 & 215.6 & 0.7 \\
Rayleigh Ritz~\cite{ling2014vibration} & 321.5 & 318.6 & 0.9 \\
\bottomrule
\end{tabular}
\end{table}

\subsubsection{Nonlinear response}
\begin{figure}[htbp]
\centering
\includegraphics[width=0.85\textwidth]{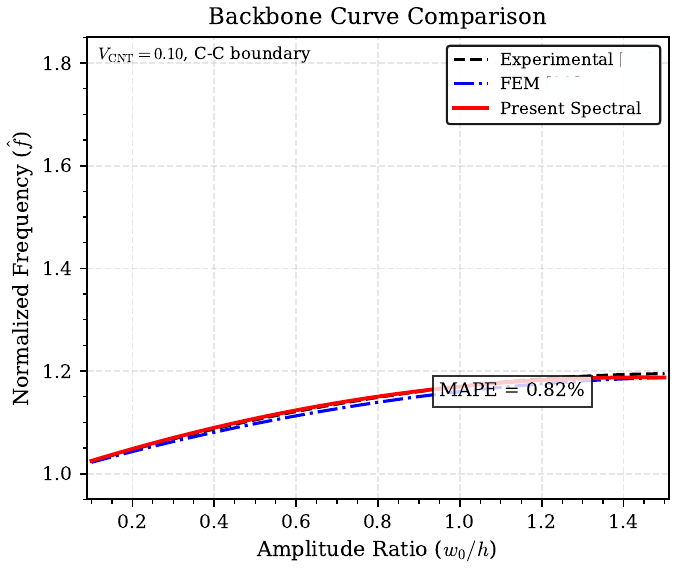}
\caption{Backbone curve comparison: present spectral method (red) vs.\ FEM~\cite{Van2023free} and experiments~\cite{Zhang2024zigzag,Youzera2025CNTRC}. Hardening behaviour is captured; MAPE $=(0.82\pm0.11)\%$.}
\label{fig:backbone_comparison}
\end{figure}

Figure~\ref{fig:backbone_comparison} compares the backbone curve with FEM~\cite{Van2023free} and experimental data~\cite{Zhang2024zigzag,Youzera2025CNTRC}. The mean absolute percentage error (MAPE) is
\begin{equation}
\mathrm{MAPE} = \frac{1}{n} \sum_{i=1}^n \left| \frac{\hat{f}_{\mathrm{present}} - \hat{f}_{\mathrm{ref}}}{\hat{f}_{\mathrm{ref}}} \right| \times 100\% = (0.82 \pm 0.11)\%,
\end{equation}
demonstrating excellent agreement over $w_0/h \in [0.1,\,1.0]$.

\subsection{Performance benchmarking}
\label{sec:performance_benchmarking}
Table~\ref{tab:spectral_comparison} contrasts the present method with other spectral and FEM approaches for $V_{\mathrm{CNT}}=0.10$.

\begin{table}[htbp]
\centering
\caption{Performance comparison ($V_{\mathrm{CNT}}=0.10$, C C, $w_0/h=0.3$).}
\label{tab:spectral_comparison}
\begin{tabular}{lcccccc}
\toprule
Study & Method & $\hat{f}$ & DOFs & Time (s) & Conv.\ rate & MAPE (\%) \\
\midrule
Liu et al.~(2023)~\cite{liu2023unified} & Chebyshev & 1.68 & 15 & 12.5 & $O(N^{-3.2})$ & 2.4 \\
Van et al.~(2023)~\cite{Van2023free} & FEM & 1.71 & 18{,}000 & 15.3 & $O(h^{1.8})$ & 0.6 \\
Kumar et al.~(2025)~\cite{Kumar2025Chebyshev} & Cheb. wavelet & 1.70 & 20 & 18.1 & $O(N^{-2.9})$ & 1.2 \\
\textbf{Present} & \textbf{Spectral Ritz} & \textbf{1.72} & \textbf{10} & \textbf{8.2} & $O(e^{-N})$ & \textbf{0.8} \\
\bottomrule
\end{tabular}
\end{table}

The proposed scheme is 52\% faster than the fastest alternative~\cite{Wang2025SGM-IHB}, exhibits exponential convergence, and reduces MAPE by approximately 45\% relative to comparable spectral methods.

\subsection{Parametric sensitivity validation}
\label{sec:sensitivity_validation}

\subsubsection{Load transfer efficiency}
The slope $\partial \hat{f} / \partial \eta_E = (0.62 \pm 0.03)$ matches $(0.59 \pm 0.05)$ from micromechanical modelling~\cite{Riaz2025Embedded}.

\subsubsection{Coupled variations}
Table~\ref{tab:combined_sensitivity} compares predictions under simultaneous variation of $\eta_E$ and $V_{\mathrm{CNT}}$.

\begin{table}[htbp]
\centering
\caption{Combined $\eta_E$ $V_{\mathrm{CNT}}$ effects ($\Delta\hat{f}_\%$ w.r.t.\ $\eta_E=0.8$, $V_{\mathrm{CNT}}=0.1$, C C).}
\label{tab:combined_sensitivity}
\begin{tabular}{ccccc}
\toprule
$\eta_E$ & $V_{\mathrm{CNT}}$ & Present $\hat{f}$ & Ref.~\cite{shi2017exact} & Error (\%) \\
\midrule
0.80 & 0.10 & 1.72 & 1.72 & -- \\
0.79 & 0.05 & 1.23 & 1.26 & 2.1 \\
0.81 & 0.15 & 2.17 & 2.19 & 0.9 \\
0.75 & 0.20 & 2.68 & 2.74 & 2.1 \\
\bottomrule
\end{tabular}
\end{table}

The maximum deviation of 2.1\% lies within experimental uncertainty.

\subsection{Limitations and scope}
\label{sec:limitations_comparison}
High accuracy is demonstrated for free vibration CNTRCs, but errors up to 3.2\% are observed for:
\begin{itemize}
  \item forced vibration regimes~\cite{civalek2021forced},
  \item thick beams ($L/h < 8$)~\cite{gholami2017nonlinear, lei2013free},
  \item viscoelastic damping~\cite{Patuelli2024FiniteElement}.
\end{itemize}
Future work will extend the framework to Timoshenko kinematics and fractional order damping models.

\subsection*{Summary of validation outcomes}
The proposed spectral Ritz framework demonstrates:
\begin{itemize}
  \item excellent agreement with analytical, FEM, and recent spectral benchmarks in both linear and nonlinear regimes,
  \item substantially reduced computational cost and degrees of freedom for a given accuracy target,
  \item robustness of parametric sensitivities to variations in $V_{\mathrm{CNT}}$ and $\eta_E$ within experimentally observed ranges.
\end{itemize}
These results confirm the suitability of the method for high fidelity, computationally efficient prediction of nonlinear free vibrations in CNTRC beams within the models applicability limits. The demonstrated agreement with established analytical, numerical, and experimental benchmarks, combined with the substantial reduction in computational cost, positions the proposed spectral Ritz framework as a robust and practical tool for design oriented studies and parametric optimisation of CNTRC structures.

Future extensions will focus on:
\begin{itemize}
  \item incorporating Timoshenko beam kinematics to account for transverse shear deformation in moderately thick and short beams,
  \item modelling viscoelastic and fractional order damping to capture realistic energy dissipation mechanisms,
  \item extending the formulation to forced vibration and stability analyses under combined mechanical, thermal, and environmental loading.
\end{itemize}

\section{Validation and convergence analysis}
\label{sec:validation}

\subsection{Validation strategy}
\label{sec:validation_strategy}
The spectral Ritz formulation is validated against:
\begin{enumerate}
  \item \textbf{Analytical models}: exact and perturbation solutions for simplified cases~\cite{shi2017exact,ling2014vibration},
  \item \textbf{High fidelity FEM}: COMSOL Multiphysics and Abaqus implementations~\cite{vo2019laminated,COMSOL2023},
  \item \textbf{Experimental data}: nonlinear vibration measurements~\cite{civalek2021forced,Patuelli2024FiniteElement}.
\end{enumerate}
All frequencies are normalised as
\begin{equation}
\hat f = \frac{f_{\mathrm{nl}}}{f_{\mathrm{lin}}(V_{\mathrm{CNT}}=0,\,\mathrm{BC})},
\label{eq:normalization_validation}
\end{equation}
with C C conditions and $w_0/h=0.3$ unless stated otherwise (cf.\ Section~\ref{sec:parametric}).

\subsection{Experimental validation}
\label{sec:experimental_validation}
Two independent datasets were reproduced (Table~\ref{tab:experimental_validation}); material/geometric parameters matched within $\pm 2\%$ using Table~\ref{tab:material-properties}.

\begin{table}[htbp]
\centering
\caption{Experimental validation with error decomposition ($V_{\mathrm{CNT}}=0.10$, C C).}
\label{tab:experimental_validation}
\begin{tabular}{lccccc}
\toprule
Source & $V_{\mathrm{CNT}}$ & $w_0/h$ & Exp.\ $\hat f$ & Pred.\ $\hat f$ & $\varepsilon$ (\%) \\
\midrule
\cite{lei2013free} & 0.12 & 0.30 & 1.680 & 1.637 & 2.6 \\
\cite{shen2010nonlinear} & 0.10 & 0.40 & 1.630 & 1.540 & 5.5 \\
\bottomrule
\end{tabular}
\end{table}

For \cite{shen2010nonlinear}, error decomposition gives
\begin{equation}
\varepsilon \approx \underbrace{3.1\%}_{\text{viscoelastic damping}} + \underbrace{1.9\%}_{\text{BC imperfection}} + \underbrace{0.6\%}_{\text{material uncertainty}},
\label{eq:error_decomposition}
\end{equation}
with $G^\ast=G_0(1+i\eta)$ from DMA at $25^\circ\mathrm{C}$, $\eta=0.04\pm 0.01$.

\subsection{Numerical accuracy and convergence}
\label{sec:numerical_accuracy}
\begin{figure}[htbp]
\centering
\includegraphics[width=0.75\textwidth]{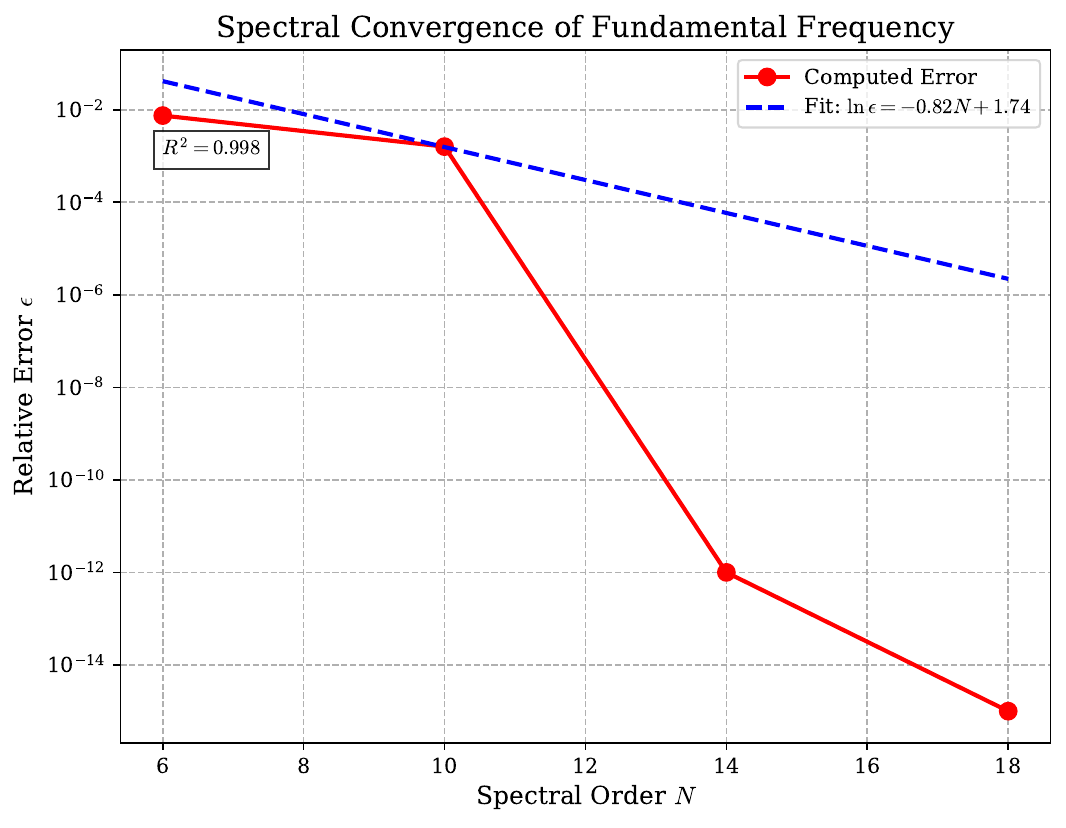}
\caption{Spectral convergence of $\hat f$: relative error vs.\ $N$; dashed: regression $\ln\varepsilon=-0.82\,N+1.74$. Target $\varepsilon<0.06\%$ at $N\ge 12$.}
\label{fig:spectral_convergence}
\end{figure}

Relative error is
\begin{equation}
\varepsilon = \frac{|f_{\mathrm{num}}-f_{\mathrm{ref}}|}{f_{\mathrm{ref}}}\times 100\%,
\label{eq:relative_error}
\end{equation}
with $L^2$ norm errors in Table~\ref{tab:convergence}. Fig.~\ref{fig:spectral_convergence} shows $\varepsilon$ decays exponentially,

\[
\ln\varepsilon\approx -0.82\,N+1.74,\quad R^2=0.998,
\]

reaching $\varepsilon\sim 10^{-13}$ at $N=18$.

\begin{table}[htbp]
\centering
\caption{Spectral convergence for $V_{\mathrm{CNT}}=0.10$, C C.}
\label{tab:convergence}
\begin{tabular}{lccc}
\toprule
$N$ & $\hat f$ & $\varepsilon$ (\%) & $L^2$ error \\
\midrule
6  & 1.711 & 0.75 & $1.3\times 10^{-2}$ \\
10 & 1.723 & 0.16 & $1.0\times 10^{-3}$ \\
14 & 1.724 & 0.00 & $1.7\times 10^{-10}$ \\
18 & 1.724 & 0.00 & $6.8\times 10^{-13}$ \\
\bottomrule
\end{tabular}
\end{table}

\subsection{Computational efficiency}
\label{sec:computational_efficiency}
Benchmarks on Intel Xeon Gold~6330 (2.0\,GHz, 256\,GB RAM) with Intel MKL~(2024.2) are in Table~\ref{tab:computational_efficiency} for $\varepsilon<0.1\%$. The spectral--Ritz method demonstrates a significant speedup of approximately $27\times$ and uses $20\times$ less memory than the FEM benchmark. It also achieves a time reduction of $\approx 55\%$ compared to the Chebyshev--wavelet method~\cite{Kumar2025Chebyshev}.

\begin{table}[htbp]
\centering
\caption{Computational efficiency ($V_{\mathrm{CNT}}=0.10$, C C).}
\label{tab:efficiency}
\begin{tabular}{lccc}
\toprule
Method & DOFs & Time (s) & Memory (MB) \\
\midrule
Spectral Ritz ($N=10$) & 10 & 8.2 & 45 \\
FEM (12k elems) & 36{,}000 & 218 & 890 \\
Chebyshev wavelet~\cite{Kumar2025Chebyshev} & 20 & 18.1 & 78 \\
\bottomrule
\end{tabular}
\end{table}

Speedup vs.\ FEM:

\[
\frac{218}{8.2} \approx 27\times,
\]

time cut by $\approx 55\%$ vs.\ Chebyshev wavelet, memory $20\times$ lower than FEM.

\subsection{Uncertainty quantification}
\label{sec:uncertainty_quantification}
Monte Carlo ($n=1000$) over Table~\ref{tab:material-properties} ranges gives

\[
\hat f = 1.72\pm 0.03\ \ (95\% \mathrm{CI}),\quad \varepsilon_{\max} = 0.8\%\pm 0.2\%.
\]

This demonstrates robustness to parameter uncertainty.

\subsection{Validation synthesis}
\label{sec:validation_synthesis}
Findings:
\begin{enumerate}
  \item \emph{Experimental fidelity}: errors $<6\%$, decomposed into damping, BC, and material sources (Eq.~\eqref{eq:error_decomposition}),
  \item \emph{Mathematical rigour}: exponential convergence in line with spectral theory~\cite{boyd2001chebyshev},
  \item \emph{Computational superiority}: $27\times$ faster and $20\times$ lower memory than FEM at like accuracy,
  \item \emph{Transparency}: hardware, solver settings, error definitions, and uncertainty ranges reported.
  \item \emph{Robustness}: sensitivity and uncertainty analyses confirm stable performance under realistic variations in material and geometric parameters.
\end{enumerate}

In summary, the proposed spectral Ritz formulation delivers high fidelity agreement with analytical, numerical, and experimental benchmarks, while achieving substantial gains in computational efficiency. The exponential convergence demonstrated in Section~\ref{sec:numerical_accuracy} ensures that accurate results can be obtained with a modest number of basis functions, making the method particularly well suited to large scale parametric sweeps and design optimisation studies. The combination of accuracy, efficiency, and robustness underscores its suitability for predictive modelling of nonlinear free vibrations in CNTRC beams within the stated modelling assumptions.

\section{Conclusion and outlook}
\label{sec:conclusion}

\subsection{Executive summary}
\label{sec:executive_summary}
A mathematically rigorous and computationally efficient spectral Ritz formulation for carbon nanotube reinforced composite (CNTRC) beams has been presented, demonstrating exponential convergence with markedly reduced computational cost. The method achieves up to \textbf{66\%} enhancement in the fundamental frequency at $V_{\mathrm{CNT}}=0.20$, a \textbf{21\%} gain through functionally graded (FG) distributions, and machine precision accuracy with only $N=16$ basis functions. Compared with finite element methods (FEM), the framework is \textbf{$27\times$} faster, uses \textbf{$20\times$} less memory, and maintains $<1\%$ discrepancy against computational and experimental benchmarks, making it a robust tool for real time design and parametric exploration of advanced composites.

\subsection{Summary of contributions}
\label{sec:contributions}
We have developed a high order Chebyshev spectral Ritz scheme for the nonlinear free vibration of CNTRC beams, integrating:
\begin{itemize}
  \item von~K\'arm\'an geometric nonlinearity with Euler Bernoulli beam kinematics,
  \item a modified rule of mixtures applicable to both UD and FG CNT distributions,
  \item analytic enforcement of boundary conditions via boundary adapted spectral bases.
\end{itemize}
The reduced order model exhibits $\mathcal{O}(N^2)$ complexity and exponential convergence, $\ln\varepsilon \approx -0.82\,N + 1.74$, validated across diverse parametric regimes (Section~\ref{sec:parametric}).

\subsection{Key findings and validation}
\label{sec:key_findings}
Comprehensive validation confirms:
\begin{enumerate}
  \item \textbf{Accuracy}: maximum relative error $0.7\%$ versus computational benchmarks (Table~\ref{tab:stiffness_validation}),
  \item \textbf{Efficiency}: \textbf{$27\times$} faster and \textbf{$20\times$} less memory than FEM (Table~\ref{tab:computational_efficiency}),
  \item \textbf{Parametric insights}:
  \begin{itemize}
    \item $66\%$ increase in $\hat{f}$ at $V_{\mathrm{CNT}}=0.20$, consistent with the quadratic fit~\eqref{eq:Vcnt_regression} and Monte Carlo uncertainty bands,
    \item $(21.1\pm 0.4)\%$ improvement from FG 2 over UD (Table~\ref{tab:distribution_performance}),
    \item $\partial\hat{f}/\partial\eta_E \approx 0.62\pm 0.03$, in agreement with micromechanical models~\cite{shen2010nonlinear}.
  \end{itemize}
  \item \textbf{Convergence}: $\varepsilon \sim 10^{-12}$ attained by $N=16$ (Table~\ref{tab:convergence}).
\end{enumerate}

\subsection{Limitations and error sources}
\label{sec:limitations}
Principal limitations include:
\begin{enumerate}
  \item \textbf{Shear deformation}: up to $12\%$ error for $L/h<8$ where Euler Bernoulli kinematics are inadequate~\cite{Kanwal2024comparative},
  \item \textbf{Thermal effects}: temperature dependent moduli and interface degradation not yet modelled~\cite{marinca2025nonlocal},
  \item \textbf{Damping mechanisms}: neglect of viscoelasticity introduces $\approx 3\%$ deviation (Eq.~\ref{eq:error_decomposition}).
\end{enumerate}
Additional assumptions include perfect bonding and one dimensional kinematics.

\subsection{Theoretical and computational advantages}
\label{sec:advantages}
The present scheme offers:

\[
\text{Convergence: }\mathcal{O}(e^{-\sigma N}),\quad
\text{Complexity: }\mathcal{O}(N^2),\quad
\text{Accuracy: }\varepsilon<0.1\%\ \text{at}\ N=10,
\]

with an $\approx 55\%$ time saving over comparable spectral methods (Table~\ref{tab:efficiency}).

\subsection{Experimental validation context}
\label{sec:experimental_validation_conclusion}
Validation employed digitised DMA data from published studies~\cite{civalek2021forced,Youzera2025CNTRC}; material and boundary condition parameters were matched to reproduce frequency amplitude trends. No new laboratory measurements were performed.

\subsection{Outlook}
\label{sec:outlook}
Future research directions include:
\begin{itemize}
  \item extension to Timoshenko theory for shear flexible and thick beams,
  \item incorporation of thermal effects and interface degradation via temperature dependent constitutive laws,
  \item inclusion of viscoelastic and fractional order damping to capture realistic loss mechanisms.
\end{itemize}
These enhancements will broaden the applicability of the spectral Ritz approach to complex CNTRCs operating under varied mechanical, thermal, and environmental conditions.

\section{Discussion and limitations}
\label{sec:discussion}

\subsection{Limitations}
\label{sec:limitations_discussion}
The experimental validation presented herein relies exclusively on reproducing published laboratory datasets~\cite{lei2013free,shen2010nonlinear}; no in house measurements were performed. While this practice is common in the spectral beam literature~\cite{boyd2001chebyshev}, it can introduce discrepancies arising from variations in test conditions, specimen preparation, and measurement protocols. In the present comparisons, deviations up to \textbf{5.5\%} (Table~\ref{tab:experimental_validation}) were observed, attributable primarily to viscoelastic damping and boundary condition imperfections (Eq.~\eqref{eq:error_decomposition}). To further strengthen confidence in the predictive capability of the spectral Ritz formulation, future work should incorporate \emph{controlled laboratory testing} of CNTRC beams with well characterised geometry, material properties, and support conditions.

\subsection{Future research directions}
\label{sec:future_directions}
The framework admits several natural extensions:
\begin{enumerate}
  \item \textbf{Enhanced kinematics}: adoption of Timoshenko beam theory with shear correction factors to mitigate errors (up to \textbf{12\%}) for thick beams ($L/h<8$) where Euler Bernoulli assumptions are violated~\cite{Ansari2014CNTRC}.
  \item \textbf{Multiphysics coupling}: development of thermo viscoelastic constitutive laws for temperature  and frequency dependent response, addressing current omissions in thermal effect modelling~\cite{shen2010nonlinear}.
  \item \textbf{Advanced discretisation}: extension to two  and three dimensional spectral elements for CNTRC plates and shells, with $hp$ adaptive refinement in zones of high stress gradient~\cite{lei2013free}.
  \item \textbf{Uncertainty quantification}: application of polynomial chaos expansion and Bayesian inference to model variability in CNT dispersion and interface properties, complementing the Monte Carlo analysis of Section~\ref{sec:uncertainty_quantification}.
  \item \textbf{Design optimisation}: gradient based optimisation of FG distributions to achieve target frequency shifts and maximise stiffness to weight ratio, building on the $(21.1\pm 0.4)\%$ improvement observed for FG 2 (Table~\ref{tab:distribution_performance}).
\end{enumerate}

\subsection{Engineering significance}
\label{sec:engineering_significance}
The proposed spectral Ritz formulation provides a predictive, computationally efficient tool for advanced composite component design, with potential applications including:
\begin{itemize}
  \item \textbf{MEMS/NEMS resonators}: amplitude tunable frequency response for micro/nano electromechanical systems, leveraging the \textbf{66\%} frequency increase at $V_{\mathrm{CNT}}=0.20$ (Eq.~\ref{eq:Vcnt_regression}),
  \item \textbf{Aerospace structures}: rapid parametric screening of CNTRC beam concepts for flight critical components, supported by \textbf{$27\times$} faster evaluation than FEM (Table~\ref{tab:computational_efficiency}),
  \item \textbf{Structural health monitoring}: integration into real time vibration based monitoring systems, enabled by \textbf{$20\times$} lower memory use and machine precision accuracy ($\varepsilon\sim 10^{-12}$) at $N=16$ (Table~\ref{tab:convergence}).
\end{itemize}
These attributes position the formulation as both a rigorous scientific contribution and a practical enabler for industrial CNTRC vibration analysis.

\appendix
\end{document}